\documentclass[12pt]{amsart}
\advance\hoffset by -0.5 truecm
\usepackage{amssymb}
\newtheorem{Theorem}{Theorem}[section]
\newtheorem{Lemma}[Theorem]{Lemma}

\newtheorem{Proposition}[Theorem]{Proposition}

\newtheorem{Remark}[Theorem]{Remark}

\def\V{\mbox{Var}}

\def\Div{\mbox{div}}
\def\Z{{\mathbb Z}}
\def\R\re
\def\V{\bf V}

\def \la{\lambda}

\def \re{{\mathbb R}}

\def \V{{\bf V}}

\def \0{\lambda_{0}}
\def \la{\lambda}
\def \ga{\gamma}

\def \G{{\mathbf G}}
\def \X{{\mathbf X}}

\newcommand{\divv}{\overset v {\operatorname{div}}\,}
\newcommand{\divh}{\overset h {\operatorname{div}}\,}
\newcommand{\divm}{\overset m {\operatorname{div}}\,}
\newcommand{\de}[2]{\frac{\partial #1}{\partial #2}}

\begin{document}
\title[Entropy production in thermostats II]{Entropy production in thermostats II}

\author[N.S. Dairbekov]{Nurlan S. Dairbekov}
\address{Laboratory of Mathematics,
Kazakh British Technical University,
Tole bi 59, 050000 Almaty, Kazakhstan }
\email{Nurlan.Dairbekov@gmail.com}

\author[G.P. Paternain]{Gabriel P. Paternain}
 \address{ Department of Pure Mathematics and Mathematical Statistics,
University of Cambridge,
Cambridge CB3 0WB, England}
 \email {g.p.paternain@dpmms.cam.ac.uk}




\begin{abstract} We show that an arbitrary Anosov Gaussian thermostat close to equilibrium has positive entropy poduction unless the external field $E$
has a global potential. The configuration space is allowed to have any
dimension and magnetic forces are also allowed.
We also show the following non-perturbative result.
Suppose a Gaussian thermostat satisfies 
\[K_{w}(\sigma)+\frac{1}{4}|E_{\sigma}|^2<0\]
for every 2-plane $\sigma$, where $K_{w}$ is the sectional curvature of the associated Weyl connection and $E_{\sigma}$ is the orthogonal projection of $E$
onto $\sigma$. Then the entropy production of any SRB measure
is positive unless $E$ has a global potential.
A related non-perturbative result is also obtained for certain generalized
thermostats on surfaces.

\end{abstract}

\maketitle

\section{Introduction}

In this paper we consider the dynamical system given by the motion
of a particle of unit mass on a closed Riemannian $n$-manifold $M$ subject
to the action of an external field $E$. We also enforce as a constraint
that the kinetic energy is a constant of motion, so the resulting equation is:
\begin{equation}
\frac{D\dot{\gamma}}{dt}=E(\ga)-\frac{\langle E(\ga),\dot{\ga}\rangle}{|\dot{\ga}|^2}\,\dot{\ga},
\label{eqt}
\end{equation}
where $D$ denotes covariant derivative and $\ga:\re\to M$ is a curve
in $M$.
This equation defines a flow $\phi$ on the unit sphere bundle $SM$ of $M$ which reduces to the geodesic flow (free motion) when $E=0$.
The kinetic energy is held fixed by Gauss' principle of least constraint and thus the system defined by (\ref{eqt}) is referred to as
{\it Gaussian thermostat}.
Thermostats have become quite popular as models in nonequilibrium statistical
mechanics \cite{CELS,Ga,GaRu,H,Ru1}. Like geodesic flows, they are
{\it reversible}, that is, the flip
$SM\ni (x,v)\mapsto (x,-v)\in SM$ conjugates $\phi_{t}$ with $\phi_{-t}$.

Let $\G_E$ be the vector field in $SM$ that generates $\phi$. An easy calculation (cf. \cite{W1}) shows that the divergence $\Div\,\G_E$ of $\G_E$ with respect to the canonical volume form $\Theta$ of $SM$ is given by
\begin{equation}\label{div-e}
\Div\, \G_E=-(n-1)\,\theta,
\end{equation}
where $\theta$ is the 1-form dual to $E$, i.e., $\theta_{x}(v)=\langle E(x),v\rangle$ and we regard $\theta$ also as a function $\theta:TM\to\re$.
We see right away that $\phi$ does not preserve the Liouville measure
(i.e. $\Theta$) unless $E=0$. But in principle, the flow may preserve other smooth measures. In fact, it is an exercise to check that $\phi$ preserves
a smooth volume form iff $\theta$ is a {\it coboundary}, that is, iff there exists a smooth function $u:SM\to\re$ which solves the {\it cohomological equation}
\begin{equation}
\G_E(u)=\theta.
\label{co}
\end{equation}
For example, suppose $\theta$ is an exact 1-form, i.e., the external field
$E$ has a global potential $U$ and write $E=-\nabla U$. Then
$\G_E(-U\circ\pi)=\theta$, where $\pi:SM\to M$ is footpoint projection: $\pi(x,v)=x$. However, in general, one does not expect to have smooth solutions
of (\ref{co}). For example, if $\theta$ is a closed, non-exact 1-form and
every homology class in $H_1(SM,\Z)$ contains a closed orbit of $\phi$, then
there is no global solution to (\ref{co})\footnote{Note that 
$\pi^*:H^{1}(M,\re)\to H^1(SM,\re)$ is an isomorphism
for $M$ different from the 2-torus, thus $\theta$ is exact iff $\pi^*\theta$
is exact.}.

In the presence of hyperbolicity, there is a close relationship
 between (\ref{co}) and the
{\it entropy production} of an SRB state $\rho$ which we now describe.
We will say that a $\phi$-invariant measure $\rho$ is an SRB measure (or state)
if $\rho$ is ergodic and
\[h_{\rho}(\phi)=\sum\,\mbox{\rm positive Lyapunov exponents},\]
where $h_{\rho}(\phi)$ is the measure theoretic entropy of $\phi$ with respect to $\rho$. The entropy production of the state $\rho$ is 
given by (cf. \cite{Ru0})
\[e_{\phi}(\rho):=-\int \Div\, \G_E\,d\rho=-\sum\,\mbox{\rm Lyapunov exponents}.\]
D. Ruelle \cite{Ru0} observed that $e_{\phi}(\rho)\geq 0$
with equality iff
\begin{equation}
h_{\rho}(\phi)=\sum\,\mbox{\rm positive Lyapunov exponents}
=-\sum\,\mbox{\rm negative Lyapunov exponents}.
\label{e=0}
\end{equation}

Suppose now that $\phi$ is an Axiom A flow (we recall the definition
in Section 3) and let $\rho$ be an SRB state. We will see in Lemma \ref{nuevo}
that if $e_{\phi}(\rho)=0$, then $\phi$ is in fact a transitive Anosov flow
and (\ref{co}) must hold. Conversely if (\ref{co}) holds, then $\phi$
preserves a smooth measure and $\phi$ is a transitive Anosov flow.
Hence $\rho$ must be the unique invariant smooth measure and consequently
(\ref{e=0}) holds which in turn implies $e_{\phi}(\rho)=0$.

Thus, for Axiom A thermostats, $e_{\phi}(\rho)=0$ iff
there exists a smooth solution of (\ref{co}).


In Section \ref{sfinal} we will explain why a transitive Anosov thermostat
is always {\it homologically full}, i.e. every homology class
in $H_1(SM,\Z)$ contains a closed orbit.
Thus, if $\theta$ is closed, but not exact (e.g. electromotive forces), then
$e_{\rho}(\phi)>0$ for any Axiom A thermostat.
This was proved by M. Wojtkowski \cite[Proposition 3.1]{W1}
assuming that $\phi$ is an Anosov flow topologically conjugate to
a geodesic flow,
and by F. Bonetto, G. Gentile and V. Mastropietro \cite{BGM} for the case
of a metric of constant negative curvature and $\theta$
a small harmonic 1-form.

The natural question now is: what happens for an arbitrary field $E$
which does not necessarily have local potentials?
In two degrees of freedom the problem was solved completely in \cite{DP2}: an
Anosov Gaussian thermostat has zero entropy production iff $E$ has a global potential.
The aim of the present paper is to provide similar results for $n$ degrees
 of freedom.

We note that the assumption that $\phi$ is uniformly hyperbolic
is known in the literature on nonequilibrium statistical
mechanics as the {\it chaotic hypothesis} of G. Gallavotti and E.G.D. Cohen:
for systems out of equilibrium, physically correct macroscopic results
will be obtained by assuming that the microscopic dynamics is uniformly
hyperbolic.
A system with $e_{\phi}(\rho)>0$ is sometimes referred to as {\it dissipative}. Dissipative Gaussian thermostats
provide a large class of examples to which one can apply the Fluctuation Theorem of Gallavotti
and Cohen \cite{GC,GC1,Ga0} (extended to Anosov flows by G. Gentile \cite{G}) and this theorem is perhaps
one of the main motivations for determining precisely which thermostats are dissipative.


In our first result we will allow magnetic forces. This
involves the addition of a {\it Lorentz force} $\bf F$ to
the right hand side of (\ref{eqt}). For each $x\in M$, ${\bf F}_x:T_xM\to T_xM$
is an antisymmetric linear map such that the 2-form 
$\langle {\bf F}_x(v),w\rangle$ is closed. We will indicate this thermostat by $\phi_{E,\bf F}$. Note that $\phi_{0,{\bf F}}$ is a magnetic flow and hence it preserves the volume form $\Theta$.
Suppose $\phi_{0,{\bf F}}$ is Anosov and $E$ is an arbitrary external field.
Then for $\varepsilon$ sufficiently small and $s\in (-\varepsilon,\varepsilon)$, the flow $\phi_{sE,{\bf F}}$ is also a transitive Anosov flow.
Moreover, the map $(-\varepsilon,\varepsilon)\ni s\mapsto e(s):=e_{\phi_{sE,{\bf F}}}(\rho_{s})$ is smooth \cite{Co,Rue1,Rue2}.
It is immediate that $e'(0)=0$ and in Section 2 we will show that
$e''(0)\geq 0$ with equality iff $E$ has a global potential.
Thus we obtain:

\medskip

\noindent {\bf Theorem A.} {\it An Anosov Gaussian thermostat close
to equilibrium has zero entropy production
if and only if the external field $E$ has a global potential. Magnetic forces are allowed at equilibrium.
}

\medskip

We now explain the non-perturbative results (which do not include magnetic forces). 
Given a 2-plane $\sigma\subset T_xM$, set:
\begin{equation}\label{k-sigma}
k(\sigma):=K(\sigma)-\Div_{\sigma}E
-|E|^2+\frac{5}{4}|E_{\sigma}|^2,
\end{equation}
where $K(\sigma)$ is the sectional curvature of the 2-plane
$\sigma$, $E_{\sigma}$ is the orthogonal projection of $E$ onto $\sigma$ and
$\Div_{\sigma}E:=\langle \nabla_\xi E, \xi\rangle+\langle \nabla_\eta E, \eta\rangle$ for any orthonormal basis $\{\xi,\eta\}$ of $\sigma$.
The expression
\[K_{w}(\sigma):=K(\sigma)-\Div_{\sigma}E
-|E|^2+|E_{\sigma}|^2,\]
is precisely the sectional curvature of the Weyl connection \cite{W2}:
\[\nabla^{w}_{X}Y=\nabla_{X}Y+\langle X,E\rangle\,Y+
\langle Y,E\rangle\,X-\langle X,Y\rangle\,E.\]
Hence
\[k(\sigma)=K_{w}(\sigma)+\frac{1}{4}|E_{\sigma}|^2.\]

\medskip

\noindent {\bf Theorem B.} {\it Let $\phi$ be a Gaussian thermostat with $k<0$
 and let $\rho$ be an SRB measure.
Then $e_{\rho}(\phi)=0$
if and only if the external field $E$ has a global potential.
}

\medskip

Like in \cite{DP2} this non-perturbative result will be established by using
Pestov type identities as in \cite{CS,DS} for geodesic flows.
A closely related result about the cohomological equation 
$\G_{E}(u)=\vartheta$, where $\vartheta$ is an {\it arbitrary} 1-form is
presented in Theorem \ref{b-1}. 

We remark that in \cite[Theorem 5.1]{W2}, M. Wojtkowski has shown that
for $n\geq 3$, the condition $k<0$
implies that $\phi$ is Anosov and for $n=2$ it suffices to assume
that $K_{w}<0$.
(In general, $K_{w}<0$ only ensures that the flow has a dominated splitting.)


Our last non-perturbative result concerns a more general class of thermostats, but it will be only for $n=2$.
In principle, nothing impedes us from considering external fields acting
on the particle which are also {\it velocity dependent}.
The way to formalize this is to say that
our external field is a {\it semibasic vector field}
$E(x,v)$, that is, a smooth map $TM \ni (x,v)\mapsto E(x,v)\in TM$ such that $E(x,v)\in T_{x}M$
for all $(x,v)\in TM$. As before the equation
\[\frac{D\dot{\gamma}}{dt}=E(\ga,\dot{\ga})-\frac{\langle E(\ga,\dot{\ga}),\dot{\ga}\rangle}{|\dot{\ga}|^2}\,\dot{\ga}.\]
defines a flow $\phi$ on the unit sphere bundle $SM$. 
These generalized thermostats
are reversible as long as $E(x,v)=E(x,-v)$.

Suppose now that $M$ is a closed oriented surface. Set
$\la(x,v):=\langle E(x,v),iv\rangle$,
where $i$ indicates rotation by $\pi/2$ according to the orientation
of the surface. The evolution of the thermostat on $SM$
can now be written as
\begin{equation}
\frac{D\dot{\gamma}}{dt}=\la(\gamma,\dot{\ga})\,i\dot{\ga}.
\label{eqgt}
\end{equation}
If $\la$ does not depend on $v$, then $\phi$ is the
magnetic flow associated with the Lorentz force ${\bf F}_{x}(v)=\la(x)iv$.
If $\la$ depends linearly on $v$, we obtain the Gaussian thermostat
(\ref{eqt}).

If we fix a Riemannian metric on $M$, its conformal class determines
a complex structure.
Given a positive integer $k$, let $\mathcal H_k$ denote the space of
holomorphic sections of the
$k$-th power of the canonical line bundle. By the Riemann-Roch theorem 
this space has complex dimension $(2k-1)(g-1)$ for $k\geq 2$
and complex dimension $g$ for $k=1$, where $g$ is the genus of $M$.
(For $k=1$ we get the holomorphic 1-forms and for $k=2$ the holomorphic
quadratic differentials.) Note that the elements in $\mathcal H_k$ can be regarded
as functions on $SM$.\footnote{Sections of the $k$-th power of the canonical
line bundle can be regarded as functions on $SM$ which transform according to the rule $f(x,e^{i\varphi}v)=e^{ik\varphi}f(x,v)$.}

Recall that $\pi:SM\to M$ is
a principal $S^{1}$-fibration and we let $V$ be the infinitesimal
generator of the action of $S^1$. If $\G$ denotes the vector field that generates the geodesic flow, 
the horizontal vector field $H$ is given by the Lie bracket $H=[V,\G]$.

\medskip

\noindent {\bf Theorem C.} {\it Let $M$ be a closed oriented surface
and consider an Anosov generalized thermostat (\ref{eqgt})
determined by $\la=\Re(q)$, where $q\in \mathcal H_{k}$.
Suppose 
\[K-H(\la)+\la^2[(k+1)^2/(2k+1)]\leq 0,\]
where $K$ is the Gaussian curvature of $M$.
Then $\phi$ has zero entropy production
if and only if $\la=0$.
}

\medskip

When $K=-1$, $k=1$, and $\la$ is  sufficiently small, the theorem is proved
in \cite{BGM} using the same perturbative methods as we will use for the
proof of Theorem A. Note that for $k$ {\it odd} the flow $\phi$
is reversible, so Theorem C provides a large class of new examples
to which the Fluctuation Theorem of Gallavotti and Cohen applies.



\bigskip


\section{Derivatives of entropy production}

\subsection{The variance} Let $\phi$ be a transitive Anosov flow on
a closed manifold
$X$. We will assume that $\phi$ is weak-mixing, i.e., the equation
$F\circ\phi_t=e^{iat}F$, $a>0$, has no continuous solutions.

Let $\mu$ be a Gibbs state associated with some H\"older continuous potential.
Given a H\"older continuous function $F:X\to\re$, the {\it variance} of $F$
with respect to $\mu$ is defined as:
\[\mbox{\rm Var}_{\mu}(F):=\lim_{T\to\infty}\frac{1}{T}\int_{X}
\left(\int_{0}^{T}(F\circ\phi_t-\overline{F})\,dt\right)^2\,d\mu,\]
where
\[\overline{F}:=\int_{X}F\,d\mu.\]
This limit exists and it appears in the central limit theorem
for hyperbolic flows \cite{Ra}.
There are other equivalent ways of expressing the variance.
Let
\[\rho_{F}(t):=\int_{X}(F\circ\phi_{t}\cdot F-\overline{F}^2)\,d\mu\]
be the auto-correlation function of $F$.
Then the variance can also be expressed as (cf. \cite[Section 4]{Po}):
\[\mbox{\rm Var}_{\mu}(F)=\int_{-\infty}^{\infty}\rho_{F}(t)\,dt
=2\int_{0}^{\infty}\rho_{F}(t)\,dt.\]
In fact the Fourier transform of $\rho_{F}$
\[\hat{\rho}^{+}_{F}(w):=\int_{0}^{\infty} e^{iwt}\rho_{F}(t)\,dt\]
defined as a distribution, has a meromorphic extension to a strip
$|\Im(w)|\leq \varepsilon$ with no pole at $w=0$ \cite{Ru-1,Po1}.
The value at $w=0$ is precisely $\mbox{\rm Var}_{\mu}(F)/2$.

\subsection{Proof of Theorem A}
We first recall the setting described in the introduction.
Consider a closed Riemannian manifold and ${\bf F}$ a Lorentz force.
Suppose $\phi_{0,{\bf F}}$ is Anosov and $E$ is an arbitrary external field.
Then, by structural stability, for $\varepsilon$ sufficiently small and $s\in (-\varepsilon,\varepsilon)$, the flow $\phi_{sE,{\bf F}}$ is also a transitive (and weak-mixing) Anosov flow.

Consider the map 
\[(-\varepsilon,\varepsilon)\ni s\mapsto e(s):
=e_{\phi_{sE,{\bf F}}}(\rho_{s}).\]
It follows from the results of G. Contreras \cite{Co} or
Ruelle \cite{Rue1,Rue2} that this map is smooth.
Indeed since $\Div\,\G_{sE,{\bf F}}=-(n-1)s\,\theta$, we have
\[e(s)=(n-1)s\int_{SM}\theta\,d\rho_{s}\]
and the results in \cite{Rue1,Rue2} assert that 
$s\mapsto \int_{SM}\theta\,d\rho_s$ is smooth. Thus
\[e'(0)=(n-1)\int_{SM}\theta\,d\mu,\]
\[e''(0)=2(n-1)\left.\frac{d}{ds}\right|_{s=0}\int_{SM}\theta\,d\rho_{s}.\]
Here $\mu:=\rho_{0}$ is the Liouville measure of $SM$.
We see right away that $e'(0)=0$ since
\[\int_{SM}\theta\,d\mu=0\]
because $\theta_{x}(v)=-\theta_{x}(-v)$.
The derivative
\[\left.\frac{d}{ds}\right|_{s=0}\int_{SM}\theta\,d\rho_{s}\]
can be computed from the results in \cite{Rue1,Rue2}.
Given a smooth function $F:SM\to\re$, the derivative
\[\left.\frac{d}{ds}\right|_{s=0}\int_{SM}F\,d\rho_{s}\]
is the limit as $\omega\to 0$ with $\Im(w)>0$ of
\begin{align*}
&\int_{0}^{\infty}e^{iwt}\int_{SM}d(F\circ\phi_{t})_{(x,v)}(Y(x,v))\,d\mu(x,v)\\
&=-\int_{0}^{\infty}e^{iwt}\int_{SM}\Div\,Y(x,v) \,F(\phi_{t}(x,v))\,d\mu(x,v)
\end{align*}
where $Y$ is such that $\G_{sE,{\bf F}}=\G_{0,{\bf F}}+s\,Y$.
Since $\Div\,(\G_{sE,{\bf F}}-\G_{0,{\bf F}})=-s\,(n-1)\theta$ we see that
\[e''(0)=2(n-1)^2\,\lim_{w\to 0}\int_{0}^{\infty}e^{iwt}\int_{SM}\theta_{x}(v) \,\theta(\phi_{t}(x,v))\,d\mu(x,v).\]
As pointed out before 
\[\omega\mapsto \int_{0}^{\infty}e^{iwt}\int_{SM}\theta_{x}(v) \,\theta(\phi_{t}(x,v))\,d\mu(x,v)\]
extends to a holomorphic function near $w=0$. We can now identify
\[2\,\lim_{w\to 0}\int_{0}^{\infty}e^{iwt}\int_{SM}\theta_{x}(v) \,\theta(\phi_{t}(x,v))\,d\mu(x,v)\]
with the variance of the function $\theta$ with respect
to the Liouville measure $\mu$.
Thus
\[e''(0)=(n-1)^2\,\mbox{\rm Var}_{\mu}(\theta).\]
The variance has the wonderful property that $\mbox{\rm Var}_{\mu}(F)\geq 0$
with equality iff $F$ is a coboundary (\cite[Section 4]{Po}).
Thus we have shown that $e'(0)=0$ and $e''(0)\geq 0$
with equality iff there exists a smooth solution $u$ to the cohomological equation
\begin{equation}
\G_{0,{\bf F}}(u)=\theta.
\label{coho}
\end{equation}
But the results in \cite[Theorem B]{DP3} give a complete understanding of the cohomological equation for Anosov magnetic flows.
Indeed, there is a solution of (\ref{coho}) iff $\theta$ is an exact form. For geodesic flows (i.e. ${\bf F}=0$) this result is proved
in \cite{DS}.

Thus, unless $E$ has a global potential, $e''(0)>0$ and therefore
$e(s)$ is strictly positive for $s\neq 0$ near zero. This shows 
Theorem A.

\subsection{Some explicit calculations of 
$\mbox{\rm Var}_{\mu}(\theta)$}
Suppose $M$ is a compact locally symmetric space of negative curvature
and suppose $\theta$ is a harmonic 1-form.
If ${\bf F}=0$, then $\mbox{\rm Var}_{\mu}(\theta)$ has been
calculated by A. Katsuda and T. Sunada in \cite[Proposition 1.3]{KS}.
They show that
\[\mbox{\rm Var}_{\mu}(\theta)=\frac{2}{h\,\mbox{\rm Vol}(M)}
\int_{M}|\theta|^2\]
where $h$ is the topological entropy of the geodesic flow of $M$.
For $n=2$ and $h=1$ ($K=-1$) we obtain
\[e''(0)=\frac{2}{\mbox{\rm Vol}(M)}\int_{M}|\theta|^2.\]
With an appropriate normalization for the $L^2$-norm
of $\theta$ we recover precisely the calculation performed
in \cite[Page 687]{BGM} to compute $e''(0)$.

It is interesting to see what happens for $n=2$ and $K=-1$ if one
adds a uniform magnetic field. Suppose we take 
${\bf F}_{x}(v)=\,iv$ and $\la\in [0,\infty)$. It is well known that
for $0\leq \la<1$, the flow $\phi_{0,\la{\bf F}}$
is Anosov and for $\la=1$ we obtain the horocycle flow.
Let $\theta$ be a harmonic 1-form and for $\la\in [0,1)$ let us try
to compute $\mbox{\rm Var}_{\la,\mu}(\theta)$, the variance
of $\theta$ with respect the flow $\phi_{0,\la{\bf F}}$
and the Liouville measure $\mu$.
A direct calculation along the lines in \cite{BGM} is possible, 
but we will take a different, more economical approach that exploits
the good properties of the variance.
It is also well known (see for example \cite{P0}) that
the flow $\phi_{0,\la{\bf F}}$ is conjugate to the geodesic flow
$\phi_{0,0}$, up to a constant time scaling by $\sqrt{1-\la^2}$.
Let $f=f_{\la}:SM\to SM$ be this conjugacy and note that it is 
immediate to check that $f_{0}$ is the identity, so $f$ is
isotopic to the identity.
Since $f$ is a conjugacy:
\[df_{(x,v)}(\X_{\la})=\G_{0,\la{\bf F}}(f(x,v)),\]
where $\X_{\la}=\sqrt{1-\la^2}\,\G_{0,0}$.
Observe that $\theta_{x}(v)=\pi^*\theta(\G_{0,\la{\bf F}})(x,v)$ 
and therefore
\begin{align*}
\theta\circ f(x,v)&=\pi^*\theta(\G_{0,\la{\bf F}})(f(x,v))=
\pi^*\theta(df_{(x,v)}\X_{\la})\\
&=f^*\pi^*\theta (\X_{\la})(x,v). 
\end{align*}
Hence
\[\mbox{\rm Var}_{\la,\mu}(\theta)=
\mbox{\rm Var}_{\X_{\la},\mu}(\theta\circ f)=
\mbox{\rm Var}_{\X_{\la},\mu}(f^*\pi^*\theta (\X_{\la})).\]
Since $\theta$ is closed, $f^*\pi^*\theta$ is a closed 1-form
in $SM$. Observe that
 $\mbox{\rm Var}_{\X_{\la},\mu}(f^*\pi^*\theta (\X_{\la}))$
only depends on the cohomology class $[f^*\pi^*\theta]$
since the variance vanishes on coboundaries.
We noted before that $f$ is isotopic to the identity, thus
\begin{align*}
\mbox{\rm Var}_{\X_{\la},\mu}(f^*\pi^*\theta (\X_{\la}))
&=\mbox{\rm Var}_{\X_{\la},\mu}(\pi^*\theta(\X_{\la}))
=\mbox{\rm Var}_{\X_{\la},\mu}(\sqrt{1-\la^2}\,\theta)\\
&=(1-\la^2)\,\mbox{\rm Var}_{\X_{\la},\mu}(\theta).
\end{align*}
From the definition of the variance it follows right away that
if $F:SM\to\re$ is any function then
\[\mbox{\rm Var}_{\X_{\la},\mu}(F)=
\sqrt{1-\la^2}\,\mbox{\rm Var}_{\mu}(F)\]
which yields
\[\mbox{\rm Var}_{\X_{\la},\mu}(f^*\pi^*\theta (\X_{\la}))=
(1-\la^2)^{3/2}\,\mbox{\rm Var}_{\mu}(\theta).\]
Summarizing
\[\mbox{\rm Var}_{\la,\mu}(\theta)=
(1-\la^2)^{3/2}\,\mbox{\rm Var}_{\mu}(\theta).\]
Thus we have obtained the following formula for the second derivative
of entropy production in the presence of a uniform magnetic field
with intensity $\la$:
\[e''_{\la}(0)=(1-\la)^{3/2}\,\frac{2}{\mbox{\rm Vol}(M)}\int_{M}|\theta|^2.\]

A completely analogous formula can be obtained for compact quotients
of complex hyperbolic space with magnetic field given by the
K\"ahler 2-form.

\section{Non-perturbative results}

All the results in the this section will be based on studying 
the cohomological equation using Pestov type identities.
The results on the cohomological equation are all collected in 
Section \ref{pesti}.

\subsection{Axiom A thermostats} A closed $\phi$-invariant set
$\Lambda$ is said to be hyperbolic if  $T(SM)$ restricted to $\Lambda$
splits as $T_{\Lambda}(SM)=\re \G_E\oplus E^{u}\oplus E^{s}$ in such a way that
there are constants $C>0$ and $0<\rho<1<\eta$ such that 
for all $t>0$ we have
\[\|d\phi_{-t}|_{E^{u}}\|\leq C\,\eta^{-t}\;\;\;\;\mbox{\rm
and}\;\;\;\|d\phi_{t}|_{E^{s}}\|\leq C\,\rho^{t}.\]
The flow is Axiom A if the nonwandering set $\Omega$ is hyperbolic and
the closed orbits are dense in $\Omega$. Recall that by the Smale spectral
decomposition, $\Omega$ is a finite union of disjoint {\it basic}
hyperbolic sets.
A hyperbolic basic set is a hyperbolic set such that:
\begin{itemize}
\item the periodic orbits of $\phi|_{\Lambda}$ are dense in $\Lambda$;
\item $\phi|_{\Lambda}$ is transitive;
\item there is an open set $U\supset\Lambda$ such that $\cap_{t\in\re}\phi_{t}(U)=\Lambda$.
\end{itemize}
The flow $\phi$ is Anosov if $SM$ is a hyperbolic set. If $\phi$ is Anosov, it is also Axiom A (but not conversely, of course). Recall that there are
examples of Anosov flows for which $\Omega$ is not the whole space \cite{FW}.

\begin{Lemma} Let $\phi$ be an Axiom A thermostat and $\rho$ an SRB state.
If $e_{\rho}(\phi)=0$, then $\phi$ is a transitive Anosov flow
and there exists a smooth solution $u$ of $\G_{E}(u)=\theta$.
\label{nuevo}
\end{Lemma}

\begin{proof} Let $\Lambda$ be the basic hyperbolic set on which
$\rho$ is supported. Since $\rho$ is an SRB measure,
\cite[Theorem 5.6]{BR} implies that $\Lambda$ is an attractor,
that is, there exists an open set $U\supset \Lambda$ such that
$\Lambda=\cap_{t\geq 0}\phi_{t}(U)$.

 As pointed out in the introduction, if $e_{\rho}(\phi)=0$,
then 
\[h_{\rho}(\phi)=\sum\,\mbox{\rm positive Lyapunov exponents}
=-\sum\,\mbox{\rm negative Lyapunov exponents}.\]
Thus $\rho$ is an SRB measure for both $\phi_{t}$ and $\phi_{-t}$.
Consequently, $\Lambda$ is an attractor for both $\phi_{t}$ and $\phi_{-t}$.
This forces $\Lambda$ to be open and since it is closed, $\Lambda=SM$ and
$\phi$ is a transitive Anosov flow.

Let $J_t^s$ and $J_{t}^u$ be the stable and unstable Jacobians of $\phi$.
If $\rho$ is an SRB measure for both $\phi_t$ and $\phi_{-t}$
then the theory of Gibbs states for transitive Anosov flows (cf. \cite[Proposition 20.3.10]{KH})
implies that $-\left.\frac{d}{dt}\right|_{t=0}\log J_t^u$ and $\left.\frac{d}{dt}\right|_{t=0}\log J_t^s$
are cohomologous (and the coboundary is the derivative along the flow of a H\"older continuous function). 
It follows that $\phi$ preserves an absolutely continuous invariant measure with positive continuous density
(and this measure would have to be $\rho$). 
An application of the smooth Liv\v sic theorem 
\cite[Corollary 2.1]{LMM}
shows that $\phi$ preserves an absolutely continuous invariant measure with positive continuous density
if and only if $\phi$ preserves a smooth volume form.
But if $\phi$ preserves a smooth volume form, then there is a smooth
solution $u$ of $\G_E(u)=\theta$ as desired.

\end{proof}

\subsection{Proof of Theorem B} We are required to prove
that if $e_{\rho}(\phi)=0$, then $E$ has a global potential.
By \cite[Theorem 5.1]{W2} the condition $k<0$ implies that $\phi$
is Anosov and hence Axiom A. By Lemma \ref{nuevo}
there is a smooth solution $u$ of $\G_{E}(u)=\theta$.
Theorem \ref{teoT} implies that $\theta$ is exact as desired.

\subsection{Proof of Theorem C}

We will need some preliminaries which can all be found in \cite{GK}.
Let $L^2(SM)$ be the space of square integrable functions with respect
to the Liouville measure of $SM$.
The space $L^{2}(SM)$ decomposes into an orthogonal direct sum of
subspaces $\sum H_{n}$, $n\in\Z$, such that on $H_{n}$, $-i\,V$ is
$n$ times the identity operator.
Consider the following first order differential operator:
\[\eta_{-}:=(\G+i\,H)/2.\]
The operator $\eta_{-}$ extends to a densely defined operator from
$H_{n}$ to $H_{n-1}$ for all $n$.
If we let $C_{n}^{\infty}(SM)=H_{n}\cap C^{\infty}(SM)$, then
$\eta_{-}:C^{\infty}_{n}\to C^{\infty}_{n -1}$ is a first order elliptic
differential operator. 
The kernel of the elliptic operator $\eta_{-}$ in $C_{k}^{\infty}(SM)$ is a finite dimensional vector space
which can be identified with $\mathcal H_{k}$.
(For all these properties see \cite{GK}.)


Now take $q\in \mathcal H_k$ and let $\la:=\Re(q)$. Then $p:=V(\la)=\Re (ikq)$.
Since $\eta_{-}q=0$, we see right away that $\G(p)+HV(p)/k=0$ and hence
Theorem \ref{kten} implies that $V(\la)$ is a coboundary iff $\la=0$. 
But $V(\la)$ is the divergence of the generalized thermostat with respect to
$\Theta$ (cf. \cite[Lemma 3.2]{DP2}), so the entropy production vanishes iff $\la=0$ as desired.

\section{Pestov identity and cohomological equation for thermostats}
\label{pesti}

\subsection{Semibasic tensor fields}

Let $\pi:TM\setminus\{0\}\to M$ be the natural projection, and let
$\beta^r_s M:=\pi^*\tau^r_s M$ denote the bundle of semibasic
tensors of degree $(r,s)$, where $\tau^r_s M$ is the bundle of
tensors of degree $(r,s)$ over $M$. Sections of the bundles
$\beta^r_sM$ are called semibasic tensor fields and the space of
all smooth sections is denoted by $C^{\infty}(\beta^r_sM)$. For
such a field $T$, the coordinate representation
$$
T=(T^{i_1\dots i_r}_{j_1\dots j_s})(x,y)
$$
holds in the domain of a standard local coordinate system $(x^i,y^i)$
on $TM\setminus\{0\}$
associated with a local coordinate system $(x^i)$ in $M$.
Under a change of a local coordinate system, the components of
a semibasic tensor field are transformed by  the same formula
as those of an ordinary tensor field on $M$.

Every ``ordinary'' tensor field on $M$ defines
a semibasic tensor field by the rule $T\mapsto T\circ\pi$, so that
the space of tensor fields on $M$ can be treated as embedded in the space
of semibasic tensor fields.

For a semibasic tensor field $(T^{i_1\dots i_r}_{j_1\dots j_s})(x,y)$,
the horizontal derivative is defined by
\begin{multline*}
T^{i_1\dots i_r}_{j_1\dots j_s|k}
=\de{}{x^k}T^{i_1\dots i_r}_{j_1\dots j_s}
-\Gamma^p_{kq}y^q\de{}{y^p}T^{i_1\dots i_r}_{j_1\dots j_s}
\\
+\sum_{m=1}^r\Gamma^{i_m}_{kp}
T^{i_1\dots i_{m-1}pi_{m+1}\dots i_r}_{j_1\dots j_s}
-\sum_{m=1}^r\Gamma^{p}_{kj_m}
T^{i_1\dots i_r}_{j_1\dots j_{m-1}pj_{m+1}\dots j_s},
\end{multline*}
and the vertical derivative by
$$
T^{i_1\dots i_r}_{j_1\dots j_s\cdot k}
=\de{}{y^k}T^{i_1\dots i_r}_{j_1\dots j_s}.
$$

The operators
$$
\nabla_{|}:C^\infty(\beta^r_sM)\to C^\infty(\beta^r_{s+1}M) 
\quad\text{and}\quad
\nabla_{\cdot}:C^\infty(\beta^r_sM)\to C^\infty(\beta^r_{s+1}M)
$$
are defined as
$$
(\nabla_{|} T)^{i_1\dots i_r}_{j_1\dots j_s k}
=\nabla_{|k} T^{i_1\dots i_r}_{j_1\dots j_{s}}
:=T^{i_1\dots i_r}_{j_1\dots j_s|k}
\quad \text{and}\quad
(\nabla_{\cdot} T)^{i_1\dots i_r}_{j_1\dots j_{s}k}
=\nabla_{\cdot k} T^{i_1\dots i_r}_{j_1\dots j_{s}}
=T^{i_1\dots i_r}_{j_1\dots j_s\cdot k}. 
$$

In \cite{PSh, Sha}, the operators $\nabla_{|}$ and $\nabla_{\cdot}$
were denoted by $\overset h \nabla$
and $\overset v \nabla$ respectively.

\subsection{Thermostats and the modified horizontal derivative}

Let $(M,g)$ be a closed connected Riemannian manifold. 
Given a vector field $E$ on $M$, define
$Y\in C^\infty(\beta^1_1 M)$ by
$$
Y_{(x,y)}(\cdot)=\frac1{|y|^2}\left(\langle y,\cdot\rangle E_x
-\langle E_x,\cdot\rangle y\right).
$$
We have seen that the equation
\begin{equation}\label{eqq}
\frac{D\dot{\gamma}}{dt}=Y_{(\gamma,\dot\gamma)}(\dot\gamma)
\end{equation}
defines the Gaussian thermostat on $SM$.

Since the flow on $SM$ defined by \eqref{eqq} 
depends only on the restriction 
of $Y$ to $SM$, we may redefine $Y$ to be
$$
Y_{(x,y)}(\cdot)=\langle y,\cdot\rangle E_x
-\langle E_x,\cdot\rangle y
$$
without changing the flow on $SM$, so that from now on
\begin{equation*}
Y^i_j(x,y)=y_jE^i(x)-E_j(x)y^i.
\end{equation*}

Given $T=(T^{i_1\dots i_r}_{j_1\dots j_s})\in C^\infty(\beta^r_s M)$,
we define the modified horizontal derivative as:
$$
T^{i_1\dots i_r}_{j_1\dots j_s: k}
=T^{i_1\dots i_r}_{j_1\dots j_s|k}
+T^{i_1\dots i_r}_{j_1\dots j_s\cdot j}Y_k^j, 
$$
so that $\nabla_{:}:C^\infty(\beta^r_sM)\to C^\infty(\beta^r_{s+1}M)$. 

For convenience, we also define $\nabla^{|}$, $\nabla^{\cdot}$,
and $\nabla^{:}$ as
$$
\nabla^{|i}=g^{ij} \nabla_{|i},\quad
\nabla^{\cdot i}=g^{ij}\nabla_{\cdot i},
\quad \nabla^{:i}=g^{ij} \nabla_{:i}. 
$$

We also set
$$
\mathbf XT^{i_1\dots i_r}_{j_1\dots j_s}
=y^kT^{i_1\dots i_r}_{j_1\dots j_s: k}.
$$
In particular, if $u\in C^\infty(TM\setminus\{0\})$, then
$$
\mathbf X u=y^i u_{:i}=y^i(u_{|i}+Y^j_iu_{\cdot j}).
$$
Note that $\mathbf X$ restricted to $SM$ coincides with $\G_E$.

It easy to see that if $\gamma$ satisfies \eqref{eqq}, then
$$
(\mathbf X T)(\gamma,\dot\gamma)
=\frac D{dt} (T(\gamma,\dot\gamma)).
$$

Straightforward caluculations give:
\begin{equation}\label{yh}
Y^i_{j|k}=y_jE^i_{,k}-E_{j,k}y^i
\end{equation}
\begin{equation}\label{yv}
Y^i_{j\cdot k}=g_{jk}E^i-E_j\delta^i_k,
\end{equation}
\begin{equation}\label{yij}
y^i_{:k}=Y^i_k=y_kE^i-E_ky^i, 
\end{equation}
where $(y^i)$ stands for the semibasic vector field $(x,y)\mapsto (y^i)$.

For $V=(V^i)\in C^\infty(\beta^1_0M)$, we set
$$
\divh V:=V^i_{|i},
\quad
\divv V:=V^i_{\cdot i},
\quad
\divm V:=V^i_{:i}.
$$

We recall the Gauss--Ostrogradski\u\i{} formulas for the horizontal
and vertical divergences \cite{Sha} (see also \cite[Section 4.2]{DP3}, 
which deals with the case of Finsler metrics). 
If $V(x,y)$ is a smooth semibasic vector field
positively homogeneous of degree $\lambda$ in $y$, 
then
\begin{align}
\label{divh}
\int_{SM}\divh V\,d\mu&=0, 
\\
\int_{SM}\divv V\,d\mu&
=\int_{SM}(\lambda+n-1)\langle V,y\rangle\,d\mu, 
\label{divv}
\end{align}
where $d\mu$ is the Liouville measure on $SM$. Whence we also have
\begin{equation}\label{divm}
\int_{SM}\divm V\, d\mu
=(\lambda+n)\int_{SM}\langle E,y \rangle \langle V,y \rangle\, d\mu
-(\lambda+1)\int_{SM}\langle V,E \rangle\, d\mu, 
\end{equation}
because
$$
\divm V=V^k_{:k}=V^k_{|k}+V^k_{\cdot i}(y_kE^i-E_ky^i)
=\divh V+\divv(\langle V,y \rangle E)-(\lambda+1)\langle V,E \rangle. 
$$

\subsection{Pestov identity}
Given a function $u:SM\to\mathbb R$, we will also denote by $u$ 
its extension to a positively homogeneous function on $TM\setminus\{0\}$
(hoping that this will not yield any confusion).

We first recall commutation formulas for horizontal and vertical
derivatives \cite{Sha} (see also \cite[Lemma 4.1]{DP3}, which deals
with the case of Finsler metrics). 
If $u\in C^\infty(TM\setminus\{0\})$, then
\begin{align}\label{vv}
u_{\cdot l\cdot k}- u_{\cdot k\cdot l}&=0,
\\
\label{hv}
u_{|l\cdot k}- u_{\cdot k|l}
&=0,
\\
\label{hh}
u_{|l|k}- u_{|k|l}
&=R^i_{lk}u_{\cdot i},
\end{align}
where $R^i_{lk}=y^jR^i_{jlk}$ and $R$ is the Riemann curvature tensor.                                               

The next lemma is an analog of \cite[Lemma 4.5]{DP3}.

\begin{Lemma}\label{com}
If $u\in C^\infty(TM\setminus\{0\})$, then
\begin{align}\label{mv}
u_{:l\cdot k}-u_{\cdot k:l}
&=(g_{lk}E^i-E_l\delta^i_k)u_{\cdot i},
\\
\label{mm}
u_{:l:k}-u_{:k:l}
&=\tilde R^i_{lk} u_{\cdot i},
\end{align}
with
\begin{equation*}
\tilde R^i_{lk}
=R^i_{kl}+(Y^i_{k|l}-Y^i_{l|k})
+(Y^j_kY^i_{l\cdot j}-Y^j_lY^i_{k\cdot j}).
\end{equation*}
\end{Lemma}

\begin{proof}
We have
$$
u_{:l\cdot k}=(u_{|l}+Y_l^iu_{\cdot i})_{\cdot k}
=u_{|l\cdot k}+Y^i_{l\cdot k}u_{\cdot i}+Y^i_{l}u_{\cdot i\cdot k},
$$
whereas
$$
u_{\cdot k:l}
=u_{\cdot k|l}+Y_l^iu_{\cdot k\cdot i}.
$$
Thus,
$$
u_{:l\cdot k}-u_{\cdot k:l}
=(u_{|l\cdot k}-u_{\cdot k|l})
+Y^i_{l\cdot k}u_{\cdot i}+Y^i_l(u_{\cdot i\cdot k}-u_{\cdot k\cdot i}).
$$

Using \eqref{vv}, \eqref{hv}, and \eqref{yv}, we come to (\ref{mv}).

Further,
\begin{multline*}
u_{:l:k}=u_{:l|k}+Y_k ^j u_{:l\cdot j}
=(u_{|l}+Y_l^ju_{\cdot j})_{|k}
+Y_k^j(u_{|l}+Y_l^s u_{\cdot s})_{\cdot j}
\\
=u_{|l|k}+Y^j_{l|k}u_{\cdot j}
+Y_l^j u_{\cdot j|k}
+Y_k^ju_{|l\cdot j}
+Y_k^jY_{l\cdot j}^s u_{\cdot s}
+Y_k^jY_l^s u_{\cdot s\cdot j}.
\end{multline*}

Therefore,
\begin{multline*}
u_{:l:k}-u_{:k:l}
=(u_{|l|k}-u_{|k|l})
+(Y^j_{l|k}-Y^j_{k|l})u_{\cdot j}
+Y^j_l (u_{\cdot j|k}-u_{|k\cdot j})
\\
+Y^j_k(u_{|l\cdot j}-u_{\cdot j|l})
+(Y^j_kY^s_{l\cdot j}-Y^j_lY^s_{k\cdot j})u_{\cdot s}
+(Y^j_kY^s_l-Y^j_lY^s_k) u_{\cdot s\cdot j}.
\end{multline*}
In view of (\ref{hv}), renaming indices and regrouping, we come to (\ref{mm}).
\end{proof}

The next lemma shows a Pestov type identity for thermostats. 

\begin{Lemma}\label{pes-n}
If\/ $u\in C^\infty(TM\setminus\{0\})$ is homogeneous of degree $0$ in $y$, then
the following holds on $SM:$
\begin{align}\label{pre-pestov}
2\langle\nabla^{:} u,\nabla^{\cdot}(\mathbf X u)\rangle
&=|\nabla^{:} u|^2
+\X(\langle\nabla^{\cdot} u,\nabla^{:} u\rangle)
-\divm((\mathbf X  u)\nabla^{\cdot} u)
+\divv((\mathbf X u)\nabla^{:} u)
\\
&-\langle \tilde {\mathbf R}_y(\nabla^{\cdot} u),\nabla^{\cdot} u\rangle
-\langle E,y\rangle \langle\nabla^{\cdot} u,\nabla^{:} u\rangle
-(n-1)(\X u)\langle E,\nabla^{\cdot}u\rangle.    \nonumber
\end{align}
\end{Lemma}

\begin{proof}
With the above notations, we can write
$$
\mathbf X u= y^i u_{:i}.
$$
Therefore,
\begin{multline}\label{start}
2\langle\nabla_{\cdot}(\mathbf X u),\nabla_{:} u\rangle
-\divv((\mathbf X u)\nabla^{:} u)
=2g^{ij}(\mathbf X u)_{\cdot i} u_{:j}
-((\mathbf X u) g^{ij} u_{:j})_{\cdot_i}
\\
=g^{ij}(\mathbf X u)_{\cdot_i} u_{:j}
-(\mathbf X u) g^{ij}u_{:j\cdot i}=I-II.
\end{multline}

We rewrite the first term on the right-hand side of (\ref{start}) as follows:
\begin{align*}
I&=g^{ij}(y^k u_{:k})_{\cdot_i} u_{:j}
=g^{ij}(u_{:i}+y^k u_{:k\cdot i}) u_{:j}
\\
&=g^{ij} u_{:i} u_{:j}
+g^{ij}y^k[u_{\cdot i:k}+(u_{:k\cdot i}-u_{\cdot i:k})] u_{:j}
\\
&=|\nabla^{:} u|^2
+y^k(g^{ij} u_{\cdot i} u_{:_j})_{:k}
-y^k g^{ij}u_{\cdot i} u_{:j:k}
+g^{ij}y^k (g_{ki}E^m-E_k\delta^m_i) u_{\cdot m} u_{:j}.
\end{align*}

Next
\[
y^k(g^{ij} u_{\cdot i} u_{:_j})_{:k}
=\X(\langle\nabla^{\cdot}u,\nabla^{:}u\rangle),
\]
\begin{align*}
y^k g^{ij}u_{\cdot i} u_{:j:k}
&=y^k g^{ij} u_{\cdot i}[u_{:k:j}+(u_{:j:k}-u_{:k:j})]
\\
&=g^{ij} u_{\cdot _i}(y^k u_{:_k})_{:j}
-g^{ij}u_{\cdot i}y^k_{:j}u_{:k}
+y^kg^{ij} u_{\cdot i} \tilde R^m_{jk} u_{\cdot m}
\\
&=\langle\nabla^{\cdot} u,\nabla^{:}(\mathbf X u)\rangle
+(\X u)\langle E,\nabla^{\cdot}u\rangle
+\langle \tilde {\mathbf R}_y(\nabla^{\cdot} u),\nabla^{\cdot} u\rangle
\end{align*}
because 
$$
g^{ij}u_{\cdot i}y^k_{:j}u_{:k}
=g^{ij}u_{\cdot i}Y^k_ju_{:k}
=y^iu_{\cdot i}E^ku_{:k}
-E^iu_{\cdot i}y^ku_{:k}
=-(\X u)\langle E,\nabla^{\cdot}u\rangle,
$$
and 
\begin{align*}
g^{ij}y^k (g_{ki}E^m-E_k\delta^m_i) u_{\cdot m} u_{:j}
&=y^jE^mu_{\cdot m} u_{:j}
-g^{ij}y^k E_k u_{\cdot i} u_{:j}\\
&=(\X u)\langle E,\nabla^{\cdot}u\rangle
-\langle E,y\rangle \langle\nabla^{\cdot} u,\nabla^{:} u\rangle.
\end{align*}

Thus,
\begin{multline}\label{start1}
I=|\nabla^{:} u|^2
+\X(\langle\nabla^{\cdot} u,\nabla^{:} u\rangle)
-\langle\nabla^{\cdot} u,\nabla^{:}(\mathbf X u)\rangle
\\
-\langle \tilde {\mathbf R}_y(\nabla^{\cdot} u),\nabla^{\cdot} u\rangle
-\langle E,y\rangle \langle\nabla^{\cdot} u,\nabla^{:} u\rangle.
\end{multline}

We rewrite the second term on the right-hand side of (\ref{start}) as
\begin{align*}
II&=(\mathbf X u) g^{ij}u_{:j\cdot i}
=(\mathbf X u) g^{ij}[u_{\cdot i:j} +(u_{:j\cdot i}-u_{\cdot i:j})]
\\
&=[(\mathbf X u) g^{ij} u_{\cdot i}]_{:j}
-(\mathbf X u)_{:j} g^{ij}u_{\cdot i}
+(\mathbf X u) g^{ij}(g_{ji}E^m-E_j\delta^m_i) u_{\cdot m}.
\end{align*}

Note that
\[
[(\mathbf X u) g^{ij} u_{\cdot i}]_{:j}
=\divm((\mathbf X  u)\nabla^{\cdot} u),
\]
that
$$
(\mathbf X u)_{:j} g^{ij}u_{\cdot i}
=\langle \nabla^{\cdot}u, \nabla^{:}(\mathbf X u)\rangle,
$$
and that
\begin{equation*}
(\mathbf X u) g^{ij}(g_{ji}E^m-E_j\delta^m_i) u_{\cdot m}
=(n-1)(\X u)\langle E,\nabla^{\cdot}u\rangle.
\end{equation*}

Thus,
\begin{equation}\label{start3}
II=\divm((\mathbf X  u)\nabla^{\cdot} u)        
-\langle \nabla^{\cdot}u, \nabla^{:}(\mathbf X u)\rangle
+(n-1)(\X u)\langle E,\nabla^{\cdot}u\rangle.
\end{equation}

Inserting (\ref{start1})--(\ref{start3}) in (\ref{start}),
we come to \eqref{pre-pestov}.
\end{proof}

Note that for the curvature term in \eqref{pre-pestov} we have, 
putting $Z=\nabla^\cdot u$: 
\begin{equation}\label{rnabla}
\langle \tilde {\mathbf R}_y(Z),Z\rangle
=\langle \mathbf R_y(Z),Z\rangle
-\langle \nabla_Z E, Z\rangle
-\langle E,Z\rangle^2. 
\end{equation}
Indeed, 
\begin{align*}
\langle \tilde {\mathbf R}_y(Z),Z\rangle
&=\left[R^i_{kl}+(Y^i_{k|l}-Y^i_{l|k})
+(Y^j_kY^i_{l\cdot j}-Y^j_lY^i_{k\cdot j})\right]y^l Z^k Z_i
\\
&=\langle \mathbf R_y(Z),Z\rangle
+(y_kE^i_{,l}-E_{k,l}y^i -y_lE^i_{,k}+E_{l,k}y^i) y^l Z^k Z_i
\\
&+[(y_kE^j-E_ky^j)(g_{lj}E^i-E_l\delta^i_j)
-(y_lE^j-E_ly^j)(g_{kj}E^i-E_{k}\delta^i_j)] y^l Z^k Z_i
\\
&=\langle \mathbf R_y(Z),Z\rangle
-E^i_{,k}Z^kZ_i
-E_kE_iZ^kZ^i, 
\end{align*}
where we used the fact that $\langle Z,y \rangle=0$
by homogeneity. 

One more useful identity is: 
\begin{equation}\label{x-n-f}
\mathbf X(\nabla^{\cdot}u)
=\nabla^{\cdot}(\mathbf Xu)
-\nabla^{:}u
-\langle E,\nabla^{\cdot}u\rangle y
+\langle E,y\rangle \nabla^{\cdot}u.
\end{equation}
Indeed,
$$
\begin{aligned}
\mathbf X (u^{\cdot i})
&=y^k(g^{ij}u_{\cdot j})_{:k}
=y^k g^{ij}(u_{:k\cdot j}
-(u_{:k\cdot j}-u_{\cdot j:k}))
\\
&=g^{ij}(y^k u_{:k})_{\cdot j}
-g^{ij}u_{:j}
-g^{ij}y^k (g_{kj}E^m-E_k\delta^m_j) u_{\cdot m},
\end{aligned}
$$
and since
\begin{equation*}
g^{ij}y^k (g_{kj}E^m-E_k\delta^m_j) u_{\cdot m}
=\langle E,\nabla^{\cdot}u\rangle y^i
-\langle E,y\rangle u^{\cdot i},
\end{equation*}
we have \eqref{x-n-f}.

\subsection{Cohomological equation}
Suppose that the cohomological equation
\begin{equation}\label{cohomological}
\mathbf G_E u=\vartheta
\end{equation}
holds with a smooth function $u$ on $SM$ and a smooth 
$1$-form $\vartheta$ on $M$.
Denoting the homogeneous extension of $u$ to $TM\setminus\{0\}$ by $u$ as before,
we get
$$
\mathbf X u(x,y)=\langle F(x),y\rangle, 
$$
where $F$ is the vector field dual to $\vartheta$ with respect to the Riemannian metric. 

Integrating \eqref{pre-pestov} against the Liouville measure $d\mu$ 
and using \eqref{divv}, \eqref{divm} yields
\begin{multline*}
2\int_{SM}\langle\nabla^{:} u,\nabla^{\cdot}(\mathbf X u)\rangle\,d\mu
=\int_{SM}\Big\{|\nabla^{:} u|^2
+\X(\langle\nabla^{\cdot} u,\nabla^{:} u\rangle)
+n(\mathbf X u)^2
\\
-\langle \tilde {\mathbf R}_y(\nabla^{\cdot} u),\nabla^{\cdot} u\rangle
-\langle E,y\rangle \langle\nabla^{\cdot} u,\nabla^{:} u\rangle
-(n-2)(\X u)\langle E,\nabla^{\cdot}u\rangle\Big\}\,d\mu.    \nonumber
\end{multline*}

Note that by \eqref{div-e}
$$
\int_{SM}\X(\langle\nabla^{\cdot} u,\nabla^{:} u\rangle)\, d\mu
=(n-1)\int_{SM}\langle E,y \rangle\langle\nabla^{\cdot} u,\nabla^{:} u\rangle\, d\mu. 
$$
Therefore, 
\begin{multline}\label{pestov}
2\int_{SM}\langle\nabla^{:} u,\nabla^{\cdot}(\mathbf X u)\rangle\,d\mu
=\int_{SM}\Big\{|\nabla^{:} u|^2
+n(\mathbf X u)^2
-\langle \tilde {\mathbf R}_y(\nabla^{\cdot} u),\nabla^{\cdot} u\rangle
\\
+(n-2)\big[\langle E,y\rangle \langle\nabla^{\cdot} u,\nabla^{:} u\rangle
-(\X u)\langle E,\nabla^{\cdot}u\rangle\big]\Big\}\,d\mu.    
\end{multline}

Using \eqref{x-n-f}, we have
\begin{multline*}
\langle E,y\rangle\langle \nabla^\cdot u,\nabla^{:}u\rangle
-(\X u)\langle E,\nabla^{\cdot}u\rangle
\\
=-\langle E,y\rangle\langle \nabla^\cdot u,\mathbf X(\nabla^{\cdot}u)\rangle
+\langle E,y\rangle\langle \nabla^\cdot u,\nabla^{\cdot}(\mathbf Xu)\rangle
+\langle E,y\rangle^2 \langle \nabla^\cdot u,\nabla^{\cdot}u\rangle
-(\X u)\langle E,\nabla^{\cdot}u\rangle
\\
=-\langle E,y\rangle\langle \nabla^\cdot u,\mathbf X(\nabla^{\cdot}u)\rangle
+\langle E,y\rangle^2|\nabla^\cdot u|^2
+\divv\left\{u\langle E,y\rangle \nabla^{\cdot}(\mathbf Xu)\right\}
-u\langle E,\nabla^{\cdot}(\mathbf Xu)\rangle
\\-u\langle E,y\rangle\,\divv \{\nabla^{\cdot}(\mathbf Xu)\}
-\divv\left\{u(\X u)E\right\}
+u\langle E,\nabla^\cdot(\X u)\rangle
\\
=-\langle E,y\rangle\langle \nabla^\cdot u,\mathbf X(\nabla^{\cdot}u)\rangle
+\langle E,y\rangle^2|\nabla^\cdot u|^2
\\
+\divv\left\{u\langle E,y\rangle \nabla^{\cdot}(\mathbf Xu)
-u(\X u)E\right\}
-u\langle E,y\rangle\,\divv \{\nabla^{\cdot}(\mathbf Xu)\}.
\end{multline*}

Plugging this in \eqref{pestov} and again using \eqref{divv}, we derive: 
\begin{multline}\label{pestov-i}
2\int_{SM}\langle\nabla^{:} u,\nabla^{\cdot}(\mathbf X u)\rangle\,d\mu
=\int_{SM}\Big\{|\nabla^{:} u|^2
+n(\mathbf X u)^2
-\langle \tilde {\mathbf R}_y(\nabla^{\cdot} u),\nabla^{\cdot} u\rangle
\\
-(n-2)\langle E,y\rangle\langle \nabla^\cdot u,\mathbf X(\nabla^{\cdot}u)\rangle
+(n-2)\langle E,y\rangle^2|\nabla^\cdot u|^2\Big\}\,d\mu, 
\end{multline}
where we used the equality
$\divv[\nabla^\cdot(\mathbf X u)]=\divv F=0$. 

Since
\begin{equation*}
\int_{SM}\Big\{|\nabla^{\cdot}(\mathbf Xu)|^2-n(\X u)^2\Big\}\, d\mu
=\int_{SM}\Big\{|F|^2-n\langle F,y\rangle^2\Big\}\,d\mu=0, 
\end{equation*}
we can rewrite \eqref{pestov-i} as follows, 
with $Z=\nabla^\cdot u$: 
\begin{multline*}
2\int_{SM}\langle\nabla^{:} u,F\rangle\, d\mu
=\int_{SM}\Big\{|\nabla^{:} u|^2
+|F|^2
-(n-2)\langle E,y\rangle\langle \mathbf X(Z),Z\rangle
\\
-\langle \tilde{\mathbf R}_y(Z),Z\rangle
+(n-2)\langle E,y\rangle^2|Z|^2\Big\}\,d\mu, 
\end{multline*}
or 
\begin{multline}\label{pestov-ii}
\int_{SM}\Big\{|F-\nabla^{:} u|^2
-(n-2)\langle E,y\rangle\langle \mathbf X(Z),Z\rangle
\\
-\langle \tilde{\mathbf R}_y(Z),Z\rangle
+(n-2)\langle E,y\rangle^2|Z|^2\Big\}\,d\mu=0. 
\end{multline}

From \eqref{x-n-f}, we obtain
\begin{equation*}
\langle \mathbf X(Z), Z \rangle
=\langle F-\nabla^{:}u , Z \rangle
+\langle E,y\rangle |Z|^2. 
\end{equation*}

Fixing any real parameter $\alpha$, we now rewrite \eqref{pestov-ii} as follows: 
\begin{multline*}
\int_{SM}\Big\{|F-\nabla^{:} u|^2
-2\alpha\langle E,y\rangle \langle F-\nabla^{:}u , Z \rangle
-(n-2-2\alpha)\langle E,y\rangle\langle \mathbf X(Z),Z\rangle
\\
-\langle\tilde{\mathbf R}_y(Z),Z\rangle
+(n-2-2\alpha)\langle E,y\rangle^2|Z|^2\Big\}\,d\mu=0, 
\end{multline*}
or
\begin{multline}\label{pestov-iii}
\int_{SM}\Big\{|F-\nabla^{:} u-\alpha\langle E,y\rangle Z|^2
-(n-2-2\alpha)\langle E,y\rangle\langle \mathbf X(Z),Z\rangle
\\
-\langle \tilde{\mathbf R}_y(Z),Z\rangle
+(n-2-2\alpha-\alpha^2)\langle E,y\rangle^2|Z|^2\Big\}\,d\mu=0. 
\end{multline}

Notice that
$$
2\langle E,y\rangle \langle\mathbf X(Z),Z\rangle 
=\mathbf X\left(\langle E,y\rangle|Z|^2\right)
-\mathbf X(\langle E,y\rangle)|Z|^2.
$$
A direct calculation gives
$$
\mathbf X(\langle E,y\rangle)
=\langle \nabla_y E,y\rangle+|E|^2-\langle E,y\rangle^2,  
$$
whence
$$
2\langle E,y\rangle \langle\mathbf X(Z),Z\rangle 
=\mathbf X\left(\langle E,y\rangle|Z|^2\right)
-(\langle \nabla_y E,y\rangle+|E|^2-\langle E,y\rangle^2)|Z|^2, 
$$
and therefore
\begin{align*}
2\int_{SM}\langle E,y&\rangle \langle\mathbf X(Z),Z\rangle \, d\mu\\
&=\int_{SM}\mathbf X\left(\langle E,y\rangle|Z|^2\right)\, d\mu
-\int_{SM}(\langle \nabla_y E,y\rangle+|E|^2-\langle E,y\rangle^2)|Z|^2\, d\mu
\\
&=\int_{SM}\Big\{n\langle E,y\rangle^2|Z|^2
-\langle \nabla_y E,y\rangle-|E|^2\Big\}\, d\mu. 
\end{align*}

Plugging this in \eqref{pestov-iii},  we get
\begin{multline*}
\int_{SM}\Big\{|F-\nabla^{:} u-\alpha\langle E,y\rangle Z|^2
+\frac{n-2-2\alpha}2\big[\langle \nabla_y E,y\rangle+|E|^2\big]
\\
-\langle \tilde{\mathbf R}_y(Z),Z\rangle 
-\left[\left(\frac{n-2-2\alpha}2\right)^2
+\left(\frac{n-2}2\right)^2\right]\langle E,y\rangle^2|Z|^2\Big\}\,d\mu=0. 
\end{multline*}

Changing $ \frac{n-2-2\alpha}2\mapsto\alpha$
and using \eqref{rnabla},  we deduce: 
\begin{multline}\label{pestov-iv}
\int_{SM}\Big\{|F-\nabla^{:} u-(n/2-1-\alpha)\langle E,y\rangle Z|^2
+\alpha\big[\langle \nabla_y E,y\rangle+|E|^2\big]
\\
-\langle \mathbf R_y(Z),Z\rangle 
+\langle \nabla_Z E, Z\rangle
+\langle E,Z\rangle^2
-\left[\alpha^2+(n/2-1)^2\right]\langle E,y\rangle^2|Z|^2\Big\}\,d\mu=0. 
\end{multline}

So,  if
$$
K(\sigma_{\xi,\eta})
-\langle \nabla_\xi E, \xi\rangle
-\alpha\langle \nabla_\eta E,\eta\rangle
-\alpha|E|^2
-\langle E,\xi\rangle^2
+\big[\alpha^2+(n/2-1)^2\big]\langle E,\eta\rangle^2
< 0
$$
for every $x\in M$ and every pair of orthogonal
unit vectors $\xi,\eta\in T_xM$, then $Z=0$. This means that
$u$ is a lift to $SM$ of a function $\varphi$ on $M$, $u(x,y)=\varphi(x)$, and 
the cohomological equation implies: $\theta=d\varphi$. 
Choosing $\alpha=1$ and putting
\begin{align*}
k_1(\sigma)&=K(\sigma)-\Div_\sigma E-|E|^2
+\big[1+(n/2-1)^2\big]|E_\sigma|^2,
\\
&=K_w(\sigma)+ \left(\frac n2-1\right)^2|E_\sigma|^2, 
\end{align*}
we arrive at the following: 

\begin{Theorem}\label{b-1}
Suppose $k_1<0$ and $\mathbf G_E(u)=\vartheta$. Then $\vartheta$ is exact. 
\end{Theorem}

It is interesting to notice that for $n=2$, $k_1=K-\Div\,E=K_w$, and that
for $n=3$,  $k_1$ equals $k$ of \eqref{k-sigma}. 

\subsection{Using the invariant measure}
The measure $f\mu$ is invariant if
$\mathbf G_E(\log f)=(n-1)\theta$. 
We let $u=\log f$,  so that
\begin{equation}\label{co-eq}
\mathbf G_E(u)=(n-1)\theta
\end{equation}
and $\nu=e^u\mu$ is a flow invariant measure. 

Let $V(x,y)$ be a smooth semibasic vector field
positively homogeneous of degree $\lambda$ in $y$. 
Since
$$
e^u\divv V
=e^u V^i_{\cdot i}=\divv(e^u V)
-\langle V,\nabla^{\cdot} u \rangle e^u, 
$$
we have by \eqref{divv}
\begin{equation}\label{divv-n}
\int_{SM}\divv V\,d\nu
=\int_{SM}\Big\{(\lambda+n-1) \langle V,y\rangle
-\langle V,\nabla^{\cdot} u \rangle\Big\}\, d\nu, 
\end{equation}
and, since
$$
e^u\divm V
=e^u V^i_{: i}=\divm(e^u V)-\langle V,\nabla^:u \rangle e^u, 
$$
we have by \eqref{divm}
\begin{equation}\label{divm-n}
\int_{SM}\divm V\,d\nu=
\int_{SM}\Big\{(\lambda+n)\langle E,y \rangle \langle V,y \rangle
-(\lambda+1)\langle V,E \rangle\, d\nu
-\langle V,\nabla^:u \rangle\Big\}\,d\nu.
\end{equation}

Let us integrate \eqref{pre-pestov} against $d\nu$. 
Using the flow invariance of $\nu$ together with \eqref{divv-n}
and \eqref{divm-n} yields: 
\begin{multline}\label{pestov-inv}
2\int_{SM}\langle\nabla^{:} u,\nabla^{\cdot}(\mathbf X u)\rangle\, d\nu
=\int_{SM}\Big\{|\nabla^{:} u|^2
+n(\mathbf X u)^2
-(n-2)(\X u)\langle E,\nabla^{\cdot}u\rangle
\\
-\langle E,y\rangle \langle\nabla^{\cdot} u,\nabla^{:} u\rangle
-\langle \tilde{\mathbf R}_y(\nabla^{\cdot} u),\nabla^{\cdot} u\rangle
\Big\}\,d\nu. 
\end{multline}
We have
\begin{multline*}
n\int_{SM}(\mathbf Xu)^2\, d\nu
=n(n-1)^2\int_{SM}\langle E,y \rangle^2 e^u\, d\mu
=(n-1)^2n\int_{SM}\langle E,y \rangle (e^u\langle E,y \rangle)\, d\mu
\\
=(n-1)^2\int_{SM}\divv(e^u\langle E,y \rangle E)\, d\mu
=(n-1)^2\int_{SM}(\langle E,y \rangle\langle \nabla^{\cdot}u,E \rangle+|E|^2)\, d\nu. 
\end{multline*}

Plugging this, \eqref{rnabla} and \eqref{co-eq} in \eqref{pestov-inv}, 
we have,  with $Z=\nabla^\cdot u$: 
\begin{multline*}
2(n-1)\int_{SM}\langle \nabla^:u,E \rangle\, d\nu
=\int_{SM}\Big\{|\nabla^:u|^2
+(n-1)\theta \langle E,Z \rangle+(n-1)^2|E|^2
\\
-\theta \langle Z,\nabla^:u\rangle
-\langle \mathbf R_y(Z),Z\rangle
+\langle \nabla_{Z}E, Z\rangle
+\langle E,Z\rangle^2
\Big\}\,d\nu=0, 
\end{multline*}
or
\begin{multline}\label{final-inv}
\int_{SM}\Big\{|\nabla^:u-(n-1)E-(1/2)\theta Z|^2
-\frac 14\theta^2|Z|^2
\\
-\langle \mathbf R_y(Z),Z\rangle
+\langle \nabla_{Z}E, Z\rangle
+\langle E,Z\rangle^2
\Big\}\,d\nu=0. 
\end{multline}

So,  if
$$
K(\sigma_{\xi,\eta})-\langle \nabla_{\xi}E, \xi\rangle
-\langle E,\xi\rangle^2+\frac 14\langle E,\eta \rangle^2<0
$$
for every $x\in M$ and every pair of orthogonal
unit vectors $\xi,\eta\in T_xM$, then $Z=0$. Passing by we note
that this condition also implies that $\phi$ is
Anosov by \cite[Theorem 4.1]{W1}.

Using \eqref{x-n-f},  we have
\begin{equation}\label{xz-e}
\mathbf X(Z)
=(n-1)E-\nabla^: u
-\langle E,Z\rangle y
+\theta Z.
\end{equation}
Then we can rewrite \eqref{final-inv} as
\begin{multline*}
\int_{SM}\Big\{|\mathbf X(Z)+\langle E,Z \rangle y-(1/2)\theta Z|^2 
-\frac 14\theta^2|Z|^2
\\
-\langle \mathbf R_y(Z),Z\rangle
+\langle \nabla_{Z}E, Z\rangle
+\langle E,Z\rangle^2
\Big\}\,d\nu=0
\end{multline*}
or
\begin{multline*}
\int_{SM}\Big\{|\mathbf X(Z)+\langle E,Z \rangle y+(1/2)\theta Z|^2 
-2\theta\langle \mathbf X(Z),Z \rangle
-\frac 14\theta^2|Z|^2
\\
-\langle \mathbf R_y(Z),Z\rangle
+\langle \nabla_{Z}E, Z\rangle
+\langle E,Z\rangle^2
\Big\}\,d\nu=0
\end{multline*}
or,  using 
$$
2\theta\langle \mathbf X(Z),Z \rangle
=\mathbf X(\theta|Z|^2)-|Z|^2\mathbf X\theta
=\mathbf X(\theta|Z|^2)-(\langle \nabla_y E,y\rangle+|E|^2-\theta^2)|Z|^2, 
$$
as
\begin{multline}\label{final-inv-1}
\int_{SM}\Big\{|\mathbf X(Z)+\langle E,Z \rangle y+(1/2)\theta Z|^2
\\
-\langle \mathbf R_y(Z),Z\rangle
+\langle \nabla_{Z}E, Z\rangle 
+\langle \nabla_{y}E, y\rangle|Z|^2
+|E|^2|Z|^2 
-\langle E,Z\rangle^2 -\frac54\theta^2|Z|^2
\Big\}\,d\nu=0. 
\end{multline}

Recall that
$$
k(\sigma)=K(\sigma)-\Div_\sigma E-|E|^2+\frac54|E_\sigma|^2. 
$$
Hence if $k(\sigma)<0$ for every $x$ and  every two-plane $\sigma\in T_xM$,
(\ref{final-inv-1}) implies $Z=0$ and hence we obtain:

\begin{Theorem} Suppose $k<0$ and $\G_{E}(u)=\theta$. Then
$\theta$ is exact.
\label{teoT}
\end{Theorem}
  
To complete these results we now show that if the thermostat
is transitive and the inequality
$k(\sigma)\le 0$ holds for all $\sigma$,  then $Z=0$. 
Indeed, in this case
$$
K(\sigma)-\Div_\sigma E-|E|^2+\frac54\langle E,\eta \rangle^2<0
$$
unless $\langle E,\xi \rangle=0$, where $\{\eta,\xi\}$ is an
orthonormal basis of $\sigma$.

Then \eqref{final-inv-1} yields
$$
\langle E,Z \rangle=0
$$
and
$$
\mathbf XZ+\langle E,Z \rangle y+(1/2)\theta Z=0. 
$$
Then
$$
\mathbf XZ=-\frac12\theta Z, 
$$
yielding
$$
\mathbf X (|Z|^2)=-\theta |Z|^2. 
$$
Assuming the set $\{Z\ne 0\}$ to be nonempty, we obtain on this set
$$
\mathbf X(\log |Z|^2)=-\theta=-(\mathbf Xu)/(n-1). 
$$
Consider this equation on a dense orbit. 
Then
$$
|Z|^2 =Ce^{-u/(n-1)}
$$
on every connected component of the intersection of 
this orbit with the set $\{Z\ne0\}$, with some nonzero constant $C$
depending on the component. Such a component is obviously 
open in this orbit. 
At the same time,  it is closed as the right hand side of the above
equality is separated from zero. This means that
the  whole orbit is in the set $\{Z\ne0\}$
and so the above holds on this orbit with the same nonzero constant, 
which means that $|Z|$ is separated from zero on a dense orbit, 
and consequently it is nonzero everywhere. We now show
that this is not possible.

Recall that $\nabla^{\cdot}u:=(u^{\cdot i})$ where
$u^{\cdot i}:=g^{ij}u_{\cdot j}$
and $u_{\cdot j}:=\frac{\partial u}{\partial y^{j}}$.
Fix $x_{0}\in M$ and consider the restriction $\tilde{u}$ of $u$
to $S_{x}M$. Since $S_{x}M$ is compact
there is $y_{0}\in S_{x}M$ for which $d_{y_{0}}\tilde{u}=0$. 
Since $u$ is homogeneous of degree zero
we must have $\nabla^{\cdot}u(x_{0},y_{0})=0$. 

Thus $Z=\nabla^{\cdot}u=0$ everywhere in $SM$.

Summarizing, we have proved:

\begin{Theorem} Let $\phi$ be a transitive Gaussian thermostat with
$k\leq 0$. Then $\phi$ preserves a smooth volume form if and only if
$E$ has a global potential.
\end{Theorem}

\subsection{Cohomological equation for generalized thermostats on surfaces}

Consider the thermostat $\phi$ determined by an arbitrary function
$\la\in C^{\infty}(SM)$ and let $\G_{\la}$ be its infinitesimal generator.

\begin{Theorem} Let $p\in C^{\infty}(SM)$ be such that
$\G(p)+HV(p)/k=0$ for some positive integer $k$, and suppose
\[K-H(\la)+\la^{2}[(k+1)^{2}/(2k+1)]\leq 0.\]
Then there exists $u\in C^{\infty}(SM)$ such that $\G_{\la}(u)=p$ if and
only if $p=0$.
\label{kten}
\end{Theorem}

\begin{proof} Note that $\G_{\la}=\G+\la V$.
We will use the following Pestov type integral identity proved in \cite[Equation
(13)]{DP2}. Given $u\in C^{\infty}(SM)$ we have:

\begin{align}\label{id1}
2\int_{SM} Hu\, V\G_{\la}u\,d\mu
&=\int_{SM}(\G_{\la}u)^2\,d\mu+\int_{SM}(Hu)^2\,d\mu\\
&-\int_{SM}(K-H(\lambda)+\lambda^2)(Vu)^2\,d\mu.\notag
\end{align}
Using that $\G(p)+HV(p)/k=0$ and that $H$ and $\G$ preserve the
Liouville measure we obtain:
\[\int_{SM}Hu\, V(p)\,d\mu=-\int_{SM}u\,HV(p)\,d\mu
=k\int_{SM}u\,\G(p)\,d\mu=-k\int_{SM}\G(u)\,p\,d\mu.\] Since
$\G(u)=p-\la V(u)$ we derive
\[\int_{SM} Hu\,
V\G_{\la}u\,d\mu=-k\int_{SM}p^{2}\,d\mu+k\int_{SM}\la\,V(u)\,p\,d\mu.\]
Combining the last equality with (\ref{id1}) yields
\[(2k+1)\int_{SM}p^{2}\,d\mu-2k\int_{SM}\la V(u)\,p\,d\mu\]
\[+\int_{SM}(Hu)^2\,d\mu
-\int_{SM}(K-H(\lambda)+\lambda^2)(Vu)^2\,d\mu=0.\]
We may rewrite this equality as:
\begin{align*}
&\int_{SM}\left(\sqrt{2k+1}\,p-\frac{k\la\,V(u)}{\sqrt{2k+1}}\right)^2\,d\mu
\\&-\int_{SM}\left(K-H(\la)+\la^{2}\frac{(k+1)^{2}}{2k+1}\right)(V(u))^2\,d\mu
+\int_{SM}(Hu)^2\,d\mu=0.
\end{align*}
Combining this equality with the hypotheses we obtain $Hu=0$.
Note that
\[K-H(\la)+\la^{2}
 \leq K-H(\la)+\la^{2}[(k+1)^{2}/(2k+1)]\leq 0.\]
Using $Hu=0$ in (\ref{id1}) we obtain $\G_{\la}(u)=p=0$.

\end{proof}

\begin{Remark}{\rm If $p(x,v)=q_{x}(v,\dots,v)$ where $q$ is a symmetric $k$-tensor, then
the condition $\G(p)+HV(p)/k=0$ is just saying that $q$ has zero divergence.
For such a $p$ and $k=1$ it suffices to assume that $\phi$ is Anosov \cite{DP2}.
It is unknown if the Anosov hypothesis is enough for $k\geq 2$. The problem is open even for
geodesic flows. We refer to \cite{SU} for partial results in this direction when $k=2$.}
\end{Remark}

\section{Final remarks and open problems}
\label{sfinal}

We begin with the following basic open problem (also raised by
Wojtkowski in \cite{W2}):

\medskip
{\it Let $\phi$ be an Anosov Gaussian thermostat on a closed $n$-manifold
with $n\geq 3$. Is it true that $\phi$ is transitive?}

\medskip

When $n=2$, a result of Ghys \cite{Ghy} ensures that $\phi$ is topologically
conjugate to the geodesic flow of a metric of constant negative curvature
and thus $\phi$ is transitive. If the weak stable and unstable bundles
of $\phi$ are transversal to the vertical subspace and $M$ supports an
Anosov geodesic flow, then a related result in \cite{Ghy} also shows that
$\phi$ is transitive. Recall that the vertical subspace $\mathcal V$
at $(x,v)\in SM$ is the kernel of $d\pi_{(x,v)}:T_{(x,v)}SM\to T_{x}M$.
Thus it is natural to ask:

\medskip
{\it Let $\phi$ be an Anosov Gaussian thermostat on a closed $n$-manifold
with $n\geq 3$. Are the weak stable and unstable bundles of $\phi$ always 
transversal to $\mathcal V$?}

\medskip

For $n=2$ this is proved in \cite{DP2} but the proof requires to know {\it apriori}
that $\phi$ is transitive, so both questions are intimately related.

We now make the following useful observation.

\begin{Proposition} Let $\phi$ be a transitive Anosov thermostat on
a closed manifold. Then $\phi$ is homologically full, that is, every homology
class in $H_{1}(SM,\Z)$ contains a closed orbit of $\phi$.
In particular $\phi$ is weak-mixing.
\label{full}
\end{Proposition}

\begin{proof} We would like to use Theorem 1 in \cite{S} which gives several
equivalent conditions for a transitive Anosov flow to be homologically full.
They all imply that $\phi$ is weak-mixing.
The one that we will use is the existence of a Gibbs state $\mu$ with zero
asymptotic cycle $\Phi_{\mu}$. As a Gibbs state we take the measure $m$
of maximal entropy. Since an Anosov thermostat is reversible via 
the flip $f(x,v)=(x,-v)$ and the measure of maximal entropy is unique,
we see that $f_{*}m=m$. 

An easy argument with the Gysin sequence of the sphere bundle $\pi:SM\to M$
shows that $\pi^*:H^{1}(M,\re)\to H^{1}(SM,\re)$ is an isomorphism, so given
$c\in H^{1}(SM,\re)$ let us write $c=[\pi^{*}\omega]$ where $\omega$
is a closed 1-form in $M$. Since $\pi^*\omega(\G_{E})=\omega$ we have:
\[\Phi_{m}(c)=\int_{SM}\omega\,dm.\]
But $f_{*}m=m$ and $\omega\circ f=-\omega$, hence $\Phi_{m}(c)=0$ for all
$c$.
\end{proof}

If the 1-form $\theta$ dual to $E$ is closed (but not exact) the results 
in \cite{WL} assert
that $\phi$ is conformally symplectic and that the weak stable and unstable bundles
are Lagrangian subspaces. In this case one can consider the action of $d\phi$ on the bundle
of Lagrangian subspaces and using the Maslov cycle, Proposition \ref{full} and arguments similar to those
in \cite[Chapter 2]{P1} and \cite{DP2}, one can show that $\mathcal V$ is transversal to
the weak bundles. (Details of this will appear elsewhere.) Of course, in this case we already know that the entropy production is positive, but
the transversality property may be of help in understanding the cohomological equation in general.

Besides transitivity and transversality of the weak bundles with $\mathcal V$, the other outstanding open problem is this:

\medskip

{\it Let $\phi$ be a transitive Anosov thermostat on a closed $n$-manifold $M$ with $n\geq 3$.
Let $\vartheta$ be a smooth 1-form on $M$. Suppose $u$ is a smooth solution of
\[\G_E(u)=\vartheta.\]
Is it true that $\vartheta$ is exact?}

\medskip

For $n=2$ this is proved in \cite{DP2} and Theorem \ref{b-1} provides an affirmative answer under certain curvature condition.

\end{document}